\documentclass[10pt]{article}
\usepackage{fancyheadings,amsbsy,amscd,amsfonts,amssymb,amstext,amsmath,latexsym,theorem,index,multicol}
\usepackage[all]{xy}

\mathchardef\tnode="020E 

\def\arc{
  \hbox{\kern -0.15em
  \vbox{\hrule width 3em height 0.6ex depth -0.5 ex}
  \kern -0.33em}}

\def\darc{
  \rlap{\lower0.2ex\arc}{\raise0.2ex\arc}}

\def\tarc{
  {\rlap{\arc}{\rlap{\lower0.4ex\arc}{\raise0.4ex\arc}}}}

\def\stroke#1{
  \kern 0.05em
  \rlap\arc{{\textstyle{#1}}\atop\phantom\arc}
  \kern -0.22em}

\def\dstroke#1{
  \kern 0.05em
  \rlap\darc{{\textstyle{#1}}\atop\phantom\darc}
  \kern -0.22em}

\def\tstroke#1{%
  \kern 0.05em
  \rlap\tarc{{\textstyle{#1}}\atop\phantom\tarc}
  \kern -0.22em}

\def\centerscript#1{
  \setbox0=\hbox{$\tnode$}
  \hbox to \wd0{\hss$\scriptstyle{#1}$\hss}}

\def\node{
  \def\super{}
  \def\sub{}
  \futurelet\next\dolabellednode}

  \let\sp=^
  \let\sb=_

  \def\dolabellednode{%
    \ifx\next\sb\let\next\getsub
    \else
      \ifx\next\sp\let\next\getsuper
      \else\let\next\donode
      \fi
    \fi
    \next}

  \def\getsub_#1{\def\sub{#1}\futurelet\next\dolabellednode}
  \def\getsuper^#1{\def\super{#1}\futurelet\next\dolabellednode}

  \def\donode{%
   \rlap{$\mathop{\phantom\tnode}\limits_{\centerscript{\sub}}^{\centerscript{\super}}$}\tnode}

\def\varcdn{
  \kern -0.03em\vbox{\kern -0.5ex  
  \hbox to \wd0{\hss\vrule width 0.04em depth 5.8ex\hss}
  \kern -0.3ex  \hbox{$\tnode$}}}

\newtheorem{theorem}{Theorem}[section]
\newtheorem{lemma}[theorem]{Lemma}
\newtheorem{proposition}[theorem]{Proposition}

{\theorembodyfont{\rmfamily}
\newtheorem{example}[theorem]{Example}
\newtheorem{definition}[theorem]{Definition}

}
\newcommand\Fqsq{\mathbb{F}_{q^2}}

\newcommand\cA{{\cal A}}
\newcommand\A{{\cal A}}

\newcommand\C{{\cal C}}
\newcommand\CC{\mathbb{C}}
\newcommand\RR{\mathbb{R}}
\newcommand\HH{\mathbb{H}}
\newcommand\G{\mathcal{G}}
\newcommand\mc[1]{{\mathcal{#1}}}

\newcommand\la{\langle}
\newcommand\ra{\rangle}
\newcommand\lb{\left\{}
\newcommand\rb{\right\}}

\newcommand\gen[1]{\langle #1 \rangle}
\newcommand\id{{\rm id}}

\newenvironment{proof}{\noindent {\bf Proof. }}{\hfill \qed \medskip}

\newcommand\ie{{i.e.}, }
\newcommand\qed{{\hfill $\Box$}}
\newcommand\pend{\qed}

\newcommand\Aut{{\rm Aut}\,}
\newcommand\Stab{{\rm Stab}\,}

\newcommand{\SO}{{\rm SO}}

\newcommand{\SL}{{\rm SL}}

\newcommand{\SU}{{\rm SU}}
\newcommand{\U}{{\rm U}}
\newcommand{\PSU}{{\rm PSU}}

\newcommand{\Sp}{{\rm Sp}}

\newcommand{\Spin}{{\rm Spin}}

\renewcommand\hat{\widehat}
\newcommand\diag{{\rm diag}}
\newcommand\B{{\cal B}}
\newcommand\typ{{\rm typ}}

\title{\bf Defining amalgams of compact Lie groups}
\author{Ralf Gramlich}

\begin{document}
\maketitle

\begin{abstract}
\noindent For $n \geq 2$ let $\Delta$ be a Dynkin diagram of rank $n$ and let $I = \{ 1, \ldots, n \}$ be the set of labels of $\Delta$. A group $G$ admits a {\em weak Phan system of type $\Delta$ over $\mathbb{C}$} if $G$ is generated by subgroups $U_i$, $i \in I$, which are central quotients of simply connected compact semisimple Lie groups of rank one, and contains subgroups $U_{i,j} = \gen{U_i,U_j}$, $i \neq j \in I$, which are central quotients of simply connected compact semisimple Lie groups of rank two such that $U_i$ and $U_j$ are rank one subgroups of $U_{i,j}$ corresponding to a choice of a maximal torus and a fundamental system of roots for $U_{i,j}$.
It is shown in this article that $G$ then is a central quotient of the simply connected compact semisimple Lie group whose complexification is the simply connected complex semisimple Lie group of type~$\Delta$. 
\end{abstract}


\section{Introduction}

In 1977 Kok-Wee Phan \cite{Phan:1977} gave a method for identifying a group $G$ as a quotient of the finite unitary group $\SU_{n+1}(q^2)$ by
finding a generating configuration of subgroups $$\SU_3(q^2) \quad \mbox{ and } \quad
\SU_2(q^2)\times \SU_2(q^2)$$ in $G$.  We begin by looking at the
configuration of subgroups in $\SU_{n+1}(q^2)$ to motivate our later
definition.
Suppose $n\ge 2$ and suppose $q$ is a prime power.  Consider
$G=\SU_{n+1}(q^2)$ acting as matrices on a Hermitian $(n+1)$-dimensional vector space over $\Fqsq$ with respect to an orthonormal basis and let $U_i\cong \SU_2(q^2)$, $i=1,2,\ldots,n$, be the
subgroups of $G$, represented as matrix groups with respect to the chosen orthonormal basis, corresponding to the $(2\times 2)$-blocks along the
main diagonal.  Let $T_i$ be the diagonal subgroup in $U_i$, which is a maximal torus of $U_i$ of size $q+1$.  When $q\ne 2$ the following hold for $1 \le i, j \le n$:
\begin{description}
\item[{\bf (P1)}] if $|i-j|>1$, then $[x,y]=1$ for all $x\in U_i$ and
$y\in U_j$;
\item[{\bf (P2)}] if $|i-j|=1$, then $\la U_i,U_j\ra$ is isomorphic to
$\SU_3(q^2)$; moreover $[x,y]=1$ for all $x\in T_i$ and $y\in T_j$; and
\item[{\bf (P3)}] the subgroups $U_i$, $1 \leq i \leq n$, generate $G$.
\end{description}
Suppose now $G$ is an arbitrary group containing a system of subgroups
$U_i\cong \SU_2(q^2)$, and suppose a maximal torus $T_i$ of size $q+1$
is chosen in each $U_i$.  If the conditions (P1)--(P3) above hold
for $G$, we will say that $G$ contains a {\em Phan system of type $A_n$ over $\Fqsq$}. Aschbacher called this configuration a generating system of type $I$ in \cite{Aschbacher:1977}. 

\medskip
In \cite{Phan:1977} Kok-Wee Phan proved the following result:

\medskip
\noindent
{\bf Phan's Theorem:} \\
{\em Let $q \geq 5$ and let $n \geq 3$. If $G$ contains a Phan system of type $A_n$ over $\Fqsq$,
then $G$ is isomorphic to a central quotient of $\SU_{n+1}(q^2)$.}

\medskip
In \cite{Phan:1977a} Phan proved similar results for finite groups corresponding to all simply laced Dynkin diagrams. For the second-generation proof of the classification of the finite simple groups \cite{Gorenstein/Lyons/Solomon:1994}, \cite{Gorenstein/Lyons/Solomon:1995}, \cite{Gorenstein/Lyons/Solomon:1998}, \cite{Gorenstein/Lyons/Solomon:1999}, \cite{Gorenstein/Lyons/Solomon:2002} the question was raised whether one could generalize and unify Phan's results. After a number of partially successful attempts by several people of reproving Phan's theorems (see, e.g., \cite{Das:1995}), the program described in \cite{Bennett/Gramlich/Hoffman/Shpectorov:2003} led to new proofs of some of Phan's old results, see \cite{Bennett/Shpectorov:2004}, \cite{Gramlich/Hoffman/Nickel/Shpectorov}, and to new unexpected Phan-type theorems, see \cite{Gramlich/Hoffman/Shpectorov:2003}, \cite{Gramlich:2004}.

\medskip
The purpose of the present article is to apply the methods from the program \cite{Bennett/Gramlich/Hoffman/Shpectorov:2003}, which have originally been developed for finite groups, to compact Lie groups, yielding a generalization of a result by Borovoi \cite{Borovoi:1984} on generators and relations in compact Lie groups. The methods and ideas used in this paper have been adopted from \cite{Bennett/Shpectorov:2004}, \cite{Gramlich:2004}, \cite{Gramlich}.

To be able to properly state the result, we have to fix the setting and to define some notions. Let $G$ be a simply connected compact semisimple Lie group of rank two, i.e., $G$ is isomorphic to $\SU_2(\mathbb{C}) \times \SU_2(\mathbb{C})$ or $\SU_3(\mathbb{C})$ or $\Spin_5(\mathbb{R}) \cong \U_2(\mathbb{H})$ or $G_{2,-14}$ by \cite{Helgason:1978}, see also 94.33 of \cite{Salzmann:1995}. Let $T$ be a maximal torus of $G$, let $\Sigma = \Sigma(G_{\mathbb{C}},T_{\mathbb{C}})$ be its root system, and let $\lb \alpha, \beta \rb$ be a fundamental system of roots of $\Sigma$, cf.\ \cite{Bourbaki:1989} or \cite{Hofmann/Morris:1998}. To the simple roots $\alpha$, $\beta$ corresponds a pair of semisimple subgroups $G_\alpha$ and $G_\beta$ of $G$ normalized by $T$ and isomorphic to $\SU_2(\mathbb{C}) \cong \Spin_3(\mathbb{R}) \cong \U_1(\mathbb{H})$, which is called a {\em standard pair} of $G$. 
If $\alpha$ and $\beta$ have different length, then the standard pair $(G_\alpha,G_\beta)$ is not conjugate to the standard pair $(G_\beta,G_\alpha)$, so, by convention, we assume that in a standard pair $(G_\alpha,G_\beta)$ the root $\alpha$ is shorter than the root $\beta$ if they have different lengths. A standard pair in a central quotient of $G$ is defined as the image of a standard pair of $G$ under the natural homomorphism. Note that the images of a standard pair in the quotient have the same isomorphism types as in $G$ modulo some central subgroups.

Moreover, for $n \geq 2$ let $\Delta$ be a Dynkin diagram of rank $n$ (see \cite{Bourbaki:2002} for a complete list) and let $I = \{ 1, \ldots, n \}$ be the set of labels of $\Delta$. A group $G$ admits a {\em weak Phan system of type $\Delta$ over $\mathbb{C}$} if $G$ is generated by subgroups $U_i$, $i \in I$, which are central quotients of simply connected compact semisimple Lie groups of rank one, and contains subgroups $U_{i,j} = \gen{U_i,U_j}$, $i \neq j \in I$, which are central quotients of simply connected compact semisimple Lie groups of rank two such that $(U_i,U_j)$ or  $(U_j,U_i)$ forms a standard pair in $U_{i,j}$. In particular the groups $U_i$ and $U_{i,j}$ have the following isomorphism types:
\begin{enumerate}
\item[{\rm (1)}] $U_i \cong \SU_2(\mathbb{C})$ or $U_i \cong \SO_3(\mathbb{R})$ for all $1 \leq i \leq n$;
\item[{\rm (2)}] $\gen{U_i,U_j} \cong \left\{
\begin{array}{ll}
(U_i\times U_j) / Z, & \mbox{in case $\node_i\hspace{3em}\node_j$, where $Z$ is a central subgroup of $U_i \times U_j$,}\\
\SU_3(\CC) \mbox{ or } \PSU_3(\CC), & \mbox{in case $\node_i\stroke{}\node_j$,} \\
\U_2(\HH) \mbox{ or } \SO_5(\RR), & \mbox{in case $\node_i\dstroke{<}\node_j$ or $\node_i\dstroke{>}\node_j$,} \\
G_{2,-14}, & \mbox{in case $\node_i\tstroke{<}\node_j$ or $\node_i\tstroke{>}\node_j$.}  
\end{array}\right.$
\end{enumerate}

\medskip
\noindent
{\bf Main Theorem.} \\
{\em Let $\Delta$ be a Dynkin diagram and let $G$ be a group admitting a weak Phan system of type $\Delta$ over $\CC$. Then $G$ is a central quotient of the simply connected compact semisimple Lie group whose complexification is the simply connected complex semisimple Lie group of type~$\Delta$. 
In particular, for irreducible Dynkin diagrams, the group $G$ is a central quotient of
\begin{itemize}
\item $\SU_{n+1}(\CC)$, if $\Delta = A_n$,
\item $\Spin_{2n+1}(\RR)$, if $\Delta = B_n$,
\item $\U_{n}(\HH)$, if $\Delta = C_n$,
\item $\Spin_{2n}(\RR)$, if $\Delta = D_n$,
\item $E_{6,-78}$, if $\Delta = E_6$,
\item $E_{7,-133}$, if $\Delta = E_7$,
\item $E_{8,-248}$, if $\Delta = E_8$,
\item $F_{4,-52}$, if $\Delta = F_4$.
\end{itemize}
}

While the theorem is true for all Dynkin diagrams, it is a tautology for Dynkin diagrams of rank at most two. In particular, the theorem does not yield an interesting characterization of the group $G_{2,-14}$. 

\medskip
This article is organized as follows.  In Section \ref{important} we remind the reader of the definition of a geometry and an amalgam and state some important lemmas. In Section \ref{exist} we recall the result by Borovoi \cite{Borovoi:1984} and give an alternative proof using geometric covering theory. In Section \ref{GHS3} we study Phan systems and Phan amalgams, indicate how to pass from one concept to the other and, moreover, prove a result on uniqueness of covers of Phan amalgams. In Section \ref{n=3}, finally, we classify the unique covers of Phan amalgams from Section \ref{GHS3} and prove the Main Theorem.

\medskip
\noindent
{\bf Acknowledgement:} The author would like to express his gratitude to Karl Heinrich Hofmann for offering a thorough overview over the area of compact Lie groups, for several insightful discussions and for guiding the author via e-mail through the library of the Institute at Oberwolfach. Thanks are also due to Linus Kramer and Karl-Hermann Neeb for help and information and additional literature. Moreover, the author would like to thank Christoph M\"uller, Christoph Wockel, Helge Gl\"ockner, Linus Kramer, and Karl-Hermann Neeb for proof-reading the paper. Finally, the author would like to point out that without the fruitful interaction in the Seminar Sophus Lie of the functional analysis group at the TU Darmstadt, the author would never have thought of applying his results from finite group theory to compact Lie groups.


\section{Geometries, amalgams and some lemmas} \label{important}

In this section we collect relevant definitions and results from incidence geometry and the theory of amalgams. See \cite{Ivanov/Shpectorov:2002} for a short introduction to the topic. A thorough introduction to incidence geometry can be found in \cite{Buekenhout:1995}.

\subsection*{Geometries}

\begin{definition}
A \textbf{pregeometry} $\mc{G}$ \textbf{over the set $I$} is a
triple $(X,*,\typ)$ consisting of a set $X$, a symmetric and
reflexive \textbf{incidence relation} $*$, and a surjective
\textbf{type function} $\typ:X\rightarrow I$, subject to the
following condition:
\begin{description}\item[(Pre)] If $x*y$ with $\typ(x)=\typ(y)$, then
$x=y$.
\end{description} The set $I$ is usually called the \textbf{type set}.
A \textbf{flag} in $X$ is a set of pairwise incident elements. The
\textbf{type} of a flag $F$ is the set $\typ(F):=\{\typ(x):x\in
F\}$. A \textbf{chamber} is a flag of type $I$.  The \textbf{rank}
of a flag $F$ is $|\typ(F)|$ and the \textbf{corank} is equal to
$|I\setminus\typ(F)|$. The cardinality of $I$ is called the {\bf rank} of $\G$. The pregeometry $\mc{G}$ is \textbf{connected} if the graph
$(X,*)$ is connected.

\medskip
A \textbf{geometry} is a pregeometry with
the additional property that
\begin{description}
\item[(Geo)] every flag is contained in a chamber.
\end{description}
\end{definition}

Let $\mc{G}=(X,*,\typ)$ be a pregeometry over $I$. An
\textbf{automorphism} of $\mc{G}$ is a permutation $\sigma$ of $X$
with $\typ(\sigma(x))=\typ(x)$, for all $x\in X$, and with
$\sigma(x)*\sigma(y)$ if and only if $x*y$, for all $x,y\in X$. A group $G$ of automorphisms of $\mc{G}$ is called {\bf flag-transitive} if for each pair $F$, $F'$ of flags of $\mc{G}$ with $\typ(F) = \typ(F')$ there exists a $g \in G$ with $g(F) = F'$. A group $G$ of automorphisms of $\mc{G}$ is called {\bf chamber-transitive} if for each pair $F$, $F'$ of flags of $\mc{G}$ with $\typ(F) = I = \typ(F')$ there exists a $g \in G$ with $g(F) = F'$. Flag-transitivity implies chamber-transitivity, for a geometry flag-transitivity and chamber-transitivity coincide, and a flag-transitive pregeometry containing a chamber automatically is a geometry, cf.\ \cite{Buekenhout:1995}.

Let $F$ be a flag of $\mc{G}$, say of type $J\subseteq I$.
Then the \textbf{residue $\mc{G}_F$ of $F$} is the pregeometry
$$(X',*_{|X'\times X'},\typ_{|I\setminus J})$$ over $I \backslash J$, with $$X':=\{x\in
X:F\cup\{x\} \mbox{ is a flag of $\mc{G}$ and } \typ(x)\notin\typ(F)\}.$$

\begin{definition}
Let $\G$ and $\hat\G$ be connected geometries over the
same type set and let $\phi:\hat\G\rightarrow\G$ be a {\bf
homomorphism} of geometries, \ie $\phi$ preserves the types and
sends incident elements to incident elements.  A surjective
homomorphism $\phi$ between connected geometries $\hat \G$ and
$\G$ is called a {\bf covering} if and only if for every nonempty
flag $\hat F$ in $\hat\G$ the map $\phi$ induces an
isomorphism between the residue of $\hat F$ in $\hat\G$ and the
residue of $F=\phi(\hat F)$ in $\G$.  Coverings of a geometry
correspond to the usual topological coverings of the flag complex.
 If
$\phi$ is an isomorphism, then the covering is said to be
\textbf{trivial}. A connected geometry $\mathcal{G}$ is called {\bf simply connected} if any covering $\widehat{\mathcal{G}} \rightarrow \mathcal{G}$ of that geometry is trivial.
\end{definition}

\begin{definition}
Let $I$ be a set, let $G$ be a group and let $(G_{i})_{i \in
I}$ be a family of subgroups of $G$. Then $(\sqcup_{i \in I} G/G_{i},*,\typ)$ with $\typ (G_{i}) = i$ and
\begin{description}
\item[(Cos)] $gG_{i} *
hG_{j}$ if and only if $gG_{i} \cap hG_{j} \neq \emptyset$
\end{description}
is a pregeometry over $I$,
the {\bf coset pregeometry of $G$} with respect to $(G_{i})_{i \in I}$. Since the type function is completely
determined by the indices, we also denote the coset pregeometry of
$G$ with respect to $(G_{i})_{i \in I}$ by
$$((G/G_{i})_{i \in I}, *).$$
The family $(G_{i})_{i \in
I}$ forms a chamber. A coset pregeometry that is a geometry is called a {\bf coset geometry}.
\end{definition}

\begin{definition}
A {\bf building geometry} is a coset geometry $((G/G_{i})_{i \in I}, *)$ where $G$ is a Chevalley group, $I$ is the set of labels of the corresponding Dynkin diagram and $(G_i)_{i \in I}$ is the collection of the maximal parabolic subgroups of $G$, cf.\ \cite{Springer:1998} or \cite{Tits:1974}. The concept of building geometries is equivalent to the concept of Tits buildings, see \cite{Brown:1989} or \cite{Buekenhout:1995}. 
\end{definition}

By Theorem IV.5.2 of \cite{Brown:1989} or by Theorem 13.32 of \cite{Tits:1974}, a building geometry of rank at least three is simply connected. In the present paper, we are interested in building geometries coming from simply connected complex semisimple Lie groups. For example, the building geometry of the group $\SL_{n+1}(\CC)$ is isomorphic to the complex projective geometry $\mathbb{P}(\mathbb{C}^{n+1})$. The building geometries of the groups $\Spin_{2n+1}(\CC)$, $\Sp_{2n}(\CC)$, $\Spin_{2n}(\CC)$ are isomorphic to the respective polar geometries, i.e., the incidence geometries of the totally isotropic subspaces of nondegenerate symmetric bilinear, respectively alternating bilinear forms of Witt index $n$ over $\CC$. 

\subsection*{Amalgams}

\begin{definition}
An {\bf amalgam} $\cA$ of groups is a 
set with a partial operation of multiplication and a collection of 
subsets $(H_i)_{i\in I}$, for some index set $I$, such that the 
following conditions hold: 
\begin{itemize}
\item[(1)] $\cA=\bigcup_{i\in I}H_i$; 
\item[(2)] the product $ab$ is defined if and only if $a,b\in H_i$ for 
some $i\in I$; 
\item[(3)] the restriction of the multiplication to each $H_i$ turns 
$H_i$ into a group; and 
\item[(4)] $H_i\cap H_j$ is a subgroup  in both $H_i$ and $H_j$ for all 
$i,j\in I$.  
\end{itemize}
\end{definition}

It follows that the groups $H_i$ share the same identity element, 
which is then the only identity element in $\cA$, and that 
$a^{-1}\in\cA$ is well-defined for every $a\in\cA$. Notice that the above definition of an amalgam of groups fits well into the general concept of an amalgam 
of groups, see \cite{Serre:2003}.

 An amalgam $\B=\bigcup_{i\in I}H_i$ is a {\bf quotient} 
of the amalgam $\A=\bigcup_{i\in I}G_i$ if there is a map $\pi$ 
from $\A$ to $\B$ such that, for each $G_i$, it restricts to a homomorphism from $G_i$ onto $H_i$. The amalgam $\A$ together with the homomorphism $\pi$ is called a {\bf cover} of the amalgam $\B$.  Two covers $(\A_1,\pi_1)$ and $(\A_2,\pi_2)$ of $\A$ are called {\bf 
equivalent} if there is an isomorphism $\phi$ of $\A_1$ onto $\A_2$, 
such that $\pi_1=\pi_2 \circ \phi$.

\begin{definition}
A group $H$ is called a {\bf completion} of an amalgam $\cA$ if there 
exists a map $\pi:\cA\rightarrow H$ such that
\begin{itemize} 
\item[(1)] for all $i\in I$ the restriction of $\pi$ to $H_i$ is a 
homomorphism of $H_i$ to $H$; and 
\item[(2)] $\pi(\cA)$ generates $H$.  
\end{itemize}
\end{definition}

Among all completions of $\cA$ there is a largest one which can be 
defined as the group having the following presentation: 
$$\mc{U}(\cA)=\la t_h\mid h\in\cA,\,t_xt_y=t_z,\mbox{ whenever $xy=z$ in
$\cA$}\ra.$$ 
Obviously, $\mc{U}(\cA)$ is a completion of $\cA$ since one can take $\pi$ 
to be the mapping $h\mapsto t_h$.  Every completion of $\cA$ is 
isomorphic to a quotient of $\mc{U}(\cA)$, and because of that $\mc{U}(\cA)$ is 
called the {\bf universal completion}. An amalgam $\mc{A}$ {\bf collapses} if $\mc{U}(\mc{A}) = 1$.

\begin{example} \label{robinson}
Consider the groups
\begin{eqnarray*}
G_1 & = & \gen{y, z \mid y^{-1} z y = z^2}, \\
G_2 & = & \gen{z, x \mid z^{-1} x z = x^2}, \\
G_3 & = & \gen{x, y \mid x^{-1} y x = y^2},
\end{eqnarray*}
 which are nontrivial and pairwise isomorphic. Let $\mc{A}$ be the amalgam given by $G_1$, $G_2$, $G_3$ and the intersections
\begin{eqnarray*}
G_1 \cap G_2 & = & \gen{z} \cong \mathbb{Z}, \\
G_1 \cap G_3 & = & \gen{y} \cong \mathbb{Z}, \\
G_2 \cap G_3 & = & \gen{x} \cong \mathbb{Z}.
\end{eqnarray*}
Then $\mc{U}(\mc{A}) = 1$ by Exercises 2.2.7 and 2.2.10 of \cite{Robinson:1982}, so $\mc{A}$ collapses.
$$
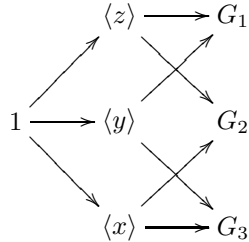
\xymatrix{
& \gen{z} \ar[r] \ar[dr] & G_1 \\
1 \ar[ur] \ar[r] \ar[dr] & \gen{y} \ar[ur] \ar[dr] & G_2 \\
& \gen{x} \ar[ur] \ar[r] & G_3
}
$$
\end{example}

\subsection*{Some lemmas}

\begin{lemma}[Tits' Lemma] \label{tits} \label{Tits2} \label{tits2}
Let $\mc{G}$ be a connected geometry over $I$ of rank at least three, let $G$ be a flag-transitive group of automorphisms of $\G$, and let $F$ be a maximal flag of $\G$.
Let $\mc{A}(\G,G,F)$ be the amalgam of stabilizers in $G$ of the elements of $F$. The geometry $\G$ is simply connected if and only if the canonical epimorphism $\mc{U}(\mc{A}(\G,G,F)) \rightarrow G$ is an isomorphism.
\end{lemma}

\begin{proof}
See Corollary 1.4.6 of \cite{Ivanov/Shpectorov:2002} or Corollary 1 of \cite{Tits:1986}.
\end{proof}

\begin{definition}
Let $\A = P_1 \cup P_2$ and $\A' = P_1' \cup P_2'$ be amalgams over an index set of cardinality two. The amalgams $\A$ and $\A'$ are of the same {\bf type} if there exist isomorphisms $\phi_i : P_i \rightarrow P_i'$ such that $\phi_i(P_1 \cap P_2) = P_1' \cap P_2'$ for $i = 1,2$.
\end{definition}

\begin{lemma}[Goldschmidt's Lemma] \label{Goldschmidt}
Let $\A = (P_1, P_2)$ be an amalgam over an index set of cardinality two, let $A_i = \Stab_{\Aut(P_i)}(P_1 \cap P_2)$ for $i = 1,2$, and let $\alpha_i : A_i \rightarrow \Aut(P_1 \cap P_2)$ be homomorphisms mapping $a \in A_i$ onto its restriction to $P_1 \cap P_2$. Then there is a one-to-one correspondence between isomorphism classes of amalgams of the same type as $\A$ and $\alpha_2(A_2)$-$\alpha_1(A_1)$ double cosets in $\Aut(P_1 \cap P_2)$. In other words, there is a one-to-one correspondence between the different isomorphism types of amalgams $P_1 \hookleftarrow (P_1 \cap P_2) \hookrightarrow P_2$ and the double cosets $\alpha_2(A_2) \backslash \Aut(P_1 \cap P_2) / \alpha_2(A_1)$.
\end{lemma}

\begin{proof}
See Lemma 2.7 of \cite{Goldschmidt:1980} or Proposition 8.3.2 of \cite{Ivanov/Shpectorov:2002}.
\end{proof}

\begin{definition}
Let $\A = (H_{i})_{i \in I}$ be an amalgam. A completion $G$ of $\A$ is called {\bf characteristic} if and only if every automorphism of $\A$ extends to an automorphism of $G$.
\end{definition}

Notice that, since $G$ is generated by the image of $\A$ under the corresponding completion map, this extension of an automorphism is unique. Clearly, the universal completion is always characteristic as is the trivial completion.

\begin{lemma}[Bennett-Shpectorov Lemma] \label{extend} \label{64}
For $i=1,2$, let $\A_i$ be an amalgam and let $G_i$ be a completion of 
$\A_i$ with completion map $\pi_i$.  Suppose there exist isomorphisms 
$\psi:\A_1\rightarrow\A_2$ and $\phi:G_1\rightarrow G_2$ such 
that $\phi \circ\pi_1=\pi_2 \circ \psi$.  If $G_1$ is a characteristic completion 
of $\A_1$, then for any isomorphism $\psi':\A_1\rightarrow\A_2$ 
there exists a unique isomorphism $\phi':G_1\rightarrow G_2$ such 
that $\phi'\circ \pi_1=\pi_2\circ \psi'$.
$$\xymatrix{
\A_1 \ar[r]^{\pi_1} \ar@/_/[d]_\psi \ar@/^/[d]^{\psi'} & G_1 \ar@/_/[d]_\phi \ar@/^/@{.>}[d]^{\phi'} \\
\A_2 \ar[r]^{\pi_2}  & G_2
}$$
\end{lemma}

\begin{proof}
See Lemma 6.4 of \cite{Bennett/Shpectorov:2004}.
\end{proof}


\section{Generators and relations} \label{exist}

Let us recall here the results by Borovoi \cite{Borovoi:1984}. Let $G$ be a simply connected compact semisimple Lie group, let $T$ be a maximal torus of $G$, let $\Sigma = \Sigma(G_{\mathbb{C}},T_{\mathbb{C}})$ be its root system, and let $\Pi$ be a system of fundamental roots of $\Sigma$. To each root $\alpha \in \Pi$ corresponds some semisimple group $G_\alpha \leq G$ of rank one such that $T$ normalizes $G_\alpha$. For simple roots $\alpha$, $\beta$, we denote by $G_{\alpha\beta}$ the group generated by the groups $G_\alpha$ and $G_\beta$, and by $\Sigma_{\alpha\beta}$ its root system relative to the torus $T_{\alpha\beta} = T \cap G_{\alpha\beta}$. The group $G_{\alpha\beta}$ is a semisimple group of rank two and $\lb \alpha, \beta \rb$ is a fundamental system of $\Sigma_{\alpha\beta}$.

Then the following assertion holds:

\begin{theorem}[Theorem of Borovoi \cite{Borovoi:1984}] \label{borovoi}
Let $G$ be a simply connected compact semisimple Lie group, let $T$ be a maximal torus of $G$, let $\Sigma = \Sigma(G_{\mathbb{C}},T_{\mathbb{C}})$ be its root system, and let $\Pi$ be a system of fundamental roots of $\Sigma$. Then the natural epimorphism $\mc{U}(\mc{A}) \rightarrow G$ is an isomorphism where $\mc{A} = (G_{\alpha\beta})_{\alpha, \beta \in \Pi}$ is the amalgam of rank one and rank two subgroups of $G$. 
\end{theorem}

Borovoi's proof consists of computations of reduced words in the group $\mc{U}(\mc{A})$ given by generators and relations. Using the theory of Tits buildings and geometric covering theory one gets the following alternative proof:

\medskip
\noindent
{\bf Geometric proof of Theorem \ref{borovoi}.} For rank at most two there is nothing to show, so we can assume that the rank is at least three. By the Iwasawa decomposition (see Theorem VI.5.1 of \cite{Helgason:1978} or Theorem III.6.32 of \cite{Hilgert/Neeb:1991}) the group $G$ acts chamber-transitively on the building geometry $\G$ of type $\Pi$ corresponding to $G_\mathbb{C}$. Let $F$ be a chamber of $\G$ stabilized by the torus $T$ of $G$, so that the stabilizers of subflags of corank one and two of $F$ with respect to the natural action of $G$ on $\G$ are exactly the groups $G_\alpha T$ and $G_{\alpha\beta}T$. By the simple connectedness of building geometries of rank at least three (cf.\ Theorem IV.5.2 of \cite{Brown:1989} or Theorem 13.32 of \cite{Tits:1974}) plus Tits' Lemma (Lemma \ref{tits}) the group $G$ equals the universal completion of the amalgam $(G_{\alpha\beta}T)_{\alpha, \beta \in \Pi}$. Finally, by Lemma 29.3 of \cite{Gorenstein/Lyons/Solomon:1995} (or by a reduction argument as in the proof of Theorem 2 of \cite{Gramlich/Hoffman/Shpectorov:2003} or in the proof of Theorem 4.3.6 of \cite{Gramlich}) the torus $T$ can be reconstructed from the rank two tori $T_{\alpha\beta}$, $\alpha, \beta \in \Pi$, and so the group $G$ actually equals the universal completion of the amalgam $(G_{\alpha\beta})_{\alpha, \beta \in \Pi}$. \pend

\begin{proposition} \label{characteristic}
Let $n \geq 2$ and let $G$ be a simply connected compact semisimple Lie group. Then the group $G$ is a characteristic completion of 
the amalgam $(G_{\alpha\beta})_{\alpha\beta \in \Pi}$ of rank one and rank two subgroups.
\end{proposition}

\begin{proof} 
By Theorem \ref{borovoi} the group $G$ is the universal completion of the amalgam $(G_{\alpha\beta})_{\alpha\beta \in \Pi}$. Therefore any automorphism of the amalgam extends to $G$, making $G$ a characteristic completion.
\end{proof}

A result similar to Theorem \ref{borovoi} has been proved by Satarov \cite{Satarov:1985} for special unitary groups over quadratic extensions of real closed fields. This case has already been covered by Borovoi's remark after his Theorem in \cite{Borovoi:1984}. Here, too, the group acts chamber-transitively on the building geometry, so our proof above applies as well.


\section{Phan systems and Phan amalgams} \label{weakPhan1}\label{GHS3}

\begin{definition}
Let $G$ be a simply connected compact semisimple Lie group of rank two, i.e., $G$ is isomorphic to $\SU_2(\mathbb{C}) \times \SU_2(\mathbb{C})$ or $\SU_3(\mathbb{C})$ or $\Spin_5(\mathbb{R}) \cong \U_2(\mathbb{H})$ or $G_{2,-14}$ by \cite{Helgason:1978}, see also 94.33 of \cite{Salzmann:1995}. Let $T$ be a maximal torus of $G$, let $\Sigma = \Sigma(G_{\mathbb{C}},T_{\mathbb{C}})$ be its root system, and let $\lb \alpha, \beta \rb$ be a fundamental root system of $\Sigma$. To the simple roots $\alpha$, $\beta$ corresponds a pair of semisimple subgroups $G_\alpha$ and $G_\beta$ of $G$ normalized by $T$ and isomorphic to $\SU_2(\mathbb{C}) \cong \Spin_3(\mathbb{R}) \cong \U_1(\mathbb{H})$ called a {\bf standard pair} of $G$. 
If $\alpha$ and $\beta$ have different length, then the standard pair $(G_\alpha,G_\beta)$ is not conjugate to the standard pair $(G_\beta,G_\alpha)$, so when speaking of a standard pair $(G_\alpha,G_\beta)$, we assume $\alpha$ to be shorter than $\beta$ if the roots have different lengths.

A standard pair in a central quotient of $G$ is defined as the image of a standard pair of $G$ under the natural homomorphism. Note that the images of a standard pair in the quotient are isomorphic to $\SU_2(\mathbb{C})$ or to $\SO_3(\RR)$.
\end{definition}

\begin{lemma} \label{7.1.1.a} \label{trans}
Standard pairs are conjugate.
\end{lemma}

\begin{proof}
This follows immediately from the fact that maximal tori are conjugate, cf.\ Theorem 6.25 of \cite{Hofmann/Morris:1998}, and the fact that, if $\alpha, \beta \in \Pi$ and $\alpha_1, \beta_1 \in \Sigma$ have the same lengths and the same angle, there exists an element $w$ of the Weyl group with $w(\alpha_1) = \alpha$ and $w(\beta_1) = \beta$, cf.\ \cite{Bourbaki:2002}.
\end{proof}

\begin{definition} \label{system}
Let $n \geq 2$, let $\Delta$ be a Dynkin diagram of rank $n$ (see \cite{Bourbaki:2002} for a complete list) and let $I = \{ 1, \ldots, n \}$ be the set of labels of $\Delta$. A group $G$ admits a {\bf weak Phan system of type $\Delta$ over $\mathbb{C}$} if $G$ is generated by subgroups $U_i \cong \SU_2(\mathbb{C})$ or $U_i \cong \SO_3(\RR)$, $i \in I$, and contains subgroups $U_{i,j} = \gen{U_i,U_j}$, $i \neq j \in I$, which are central quotients of simply connected compact semisimple Lie groups of rank two such that $(U_i,U_j)$ or  $(U_j,U_i)$ forms a standard pair in $U_{i,j}$. In particular, any $U_{i,j}$ is isomorphic to a central quotient of $\SU_2(\mathbb{C}) \times \SU_2(\mathbb{C})$ or to $\SU_3(\mathbb{C})$ or $\PSU_3(\mathbb{C})$ or $\U_2(\mathbb{H}) \cong \Spin_5(\mathbb{R})$ or $\SO_5(\mathbb{R})$ or $G_{2,-14}$ depending on the subdiagram of $\Delta$ induced on $i$ and $j$.
\end{definition}

The paramount examples for groups with a weak Phan system are the simply connected compact semisimple Lie groups together with the amalgam $(G_{\alpha\beta})_{\alpha\beta \in \Pi}$ of rank one and rank two subgroups. Any central quotient of such a group of rank at least two also admits a weak Phan system.

\begin{definition} \label{arbitraryan} \label{arbitrarycn} \label{arbitrarydn}
 A {\bf Phan amalgam} is an 
amalgam $\A=(L_{\alpha\beta})_{\alpha, \beta \in \Pi}$, where $L_{\alpha\beta}$ is a group isomorphic to a central
quotient of $G_{\alpha\beta}$ where it is required that $L_{\alpha}$ and $L_{\beta}$ are the images of $G_{\alpha}$, respectively $G_{\beta}$
under the natural epimorphism from $G_{\alpha\beta}$ onto $L_{\alpha\beta}$.
A Phan amalgam is called {\bf irreducible} if it is obtained from the natural amalgam $(G_{\alpha\beta})_{\alpha, \beta \in \Pi}$ of a simply connected compact almost simple Lie group, i.e., if the Dynkin diagram of that group is connected or, equivalently, if the corresponding root system is irreducible, cf.\ \cite{Bourbaki:2002}. A complete list of the compact almost simple Lie groups can be found in \cite{Helgason:1978} or \cite{Salzmann:1995}.
A  Phan amal\-gam is called {\bf strongly noncollapsing} if there exists a completion $\pi : \mc{A} \rightarrow G$ such that the kernel of the restriction $\pi_{|L_{\alpha_i}}$ is central for each $i \in I$.
The {\bf rank} of a Phan amalgam is defined to be the rank of the corresponding fundamental system $\Pi$.
The amalgam $(G_{\alpha\beta})_{\alpha, \beta \in \Pi}$ is called a {\bf standard Phan amalgam}. 
\end{definition}

If a group $G$ contains a weak Phan system $U_1,\dots, U_n$, then 
$\A=(U_{i,j})_{i,j\in I}$ is a strongly noncollapsing Phan amalgam.
The converse is also true: a Phan amalgam admitting a faithful completion $G$ turns the group $G$ into a group with a weak Phan system of the respective type.

\begin{definition}
A Phan amalgam $(L_{\alpha\beta})_{\alpha, \beta \in \Pi}$ is called {\bf unambiguous} if
 every $L_{\alpha\beta}$ is isomorphic to the corresponding $G_{\alpha\beta}$.
\end{definition}

\begin{proposition} \label{refinement} \label{similar}
Every Phan amalgam $\A$ has an unambiguous covering $\hat\A$ that is unique up to equivalence of coverings. 
Furthermore, every (strongly) noncollapsing Phan amalgam $\A$ has a unique (up to equivalence of coverings)
unambiguous (strongly) noncollapsing covering $\hat\A$.
\end{proposition}

\begin{proof} 
We will proceed by induction on $|S|$, where $S$ is a subset of ${\Pi \choose 1} \cup {\Pi \choose 2}$ which is closed under taking subsets and $\A=(L_J)_{J \in S}$.  Our basis is the case $S=\emptyset$ which vacuously yields an unambiguous amalgam.  Suppose now that $S$ is non-empty, and that for 
every subset $S'\subsetneq S$ the claim holds.  Let $J$ be an element of $S$ which is maximal with respect to inclusion and define $S'=S\setminus\{J\}$ and 
$\A'=(L_{J'})_{J'\in S'}$.  Then $S'$ is closed under taking subsets, and $\A'$ is a 
 subamalgam in $\A$.

By the inductive assumption, there is a unique unambiguous covering 
 amalgam $(\hat\A'=(\hat L_{J'})_{J'\in S'},\pi')$ of $\A'$.  We 
will find an unambiguous covering $(\hat\A,\pi)$ of $\A$ by gluing a 
copy of $G_J$ to $\hat\A'$ and by extending $\pi'$ to the new member 
of the amalgam.  To glue $G_J$ to the amalgam $\hat\A'$, we need to 
construct an isomorphism from the subamalgam $\hat{\mc{L}}=(\hat L_{J'})_{J'\subsetneq J}$ of $\hat\A'$ onto the corresponding amalgam 
$\mc{G}=(G_{J'})_{J'\subsetneq J}$ of subgroups of $G_J$.  By the 
definition of a Phan amalgam, there is a homomorphism $\psi$ from 
$G_J$ onto $L_J$ mapping $\mc{G}$ onto $\mc{L}=(L_{J'})_{J'\subsetneq J}$.  
 Note that $(\hat{\mc{L}},\pi'|_{\hat{\mc{L}}})$ and $(\mc{G},\psi|_{\mc{G}})$ are two 
unambiguous coverings of $\mc{L}$.  By induction, the uniqueness of the 
unambiguous covering holds so that there is an amalgam isomorphism 
$\phi$ from $\hat{\mc{L}}$ onto $\mc{G}$ such that $\psi \circ \phi =\pi'|_{\hat{\mc{L}}}$.  Clearly, 
$\phi$ tells us how to glue $G_J$ to $\hat\A'$ to produce 
$\hat\A$ and, furthermore, as $\pi$ we can take the union of $\psi$ 
and $\pi'$.  The condition $\psi \circ \phi =\pi'|_{\hat{\mc{L}}}$ guarantees that 
$\psi$ and $\pi'$ agree on the intersection $\hat{\mc{L}} \stackrel{\phi}{\cong}\mc{G}$.  Finally, notice that $\hat\A$ is an unambiguous Phan 
amalgam, so $(\hat\A,\pi)$ is an unambiguous covering of 
$\A$.
This completes the proof of the existence of an unambiguous covering 
$\hat\A$.  

Now we will prove the uniqueness.  Suppose we have two 
such coverings $\hat\B=(B_J)_{J\in S}$ and $\hat\C=(C_J)_{J\in 
S}$ with corresponding amalgam homomorphism $\pi_1$ and $\pi_2$ 
onto $\A$.  Select $J$ as an element of $S$ which is maximal with respect to inclusion, and define 
$S'=S\setminus\{J\}$.  Let $\A'$, $\hat\B'$ and $\hat\C'$ be the 
subamalgams of shape $S'$ in $\A$, $\hat\B$ and $\hat\C$, 
respectively.  By induction, there exists an isomorphism $\phi$ from 
$\hat\B'$ onto $\hat\C'$ such that $\pi_1|_{\hat\B'}=\pi_2 \circ \phi$.  
It suffices to extend $\phi$ to $B_J$.

We have to deal with two cases: First, let us assume that $J = \lb \alpha, \beta \rb$ where $\alpha$ and $\beta$ are orthogonal roots.  In this 
case, $B_{\alpha\beta} \cong C_{\alpha\beta} \cong G_{\alpha\beta}$ is isomorphic to a direct product of $B_{\alpha} \cong C_{\alpha} \cong G_{\alpha}$ and 
$B_{\beta} \cong C_{\beta} \cong G_{\beta}$.  Clearly 
$\phi$ is already known on $B_{\alpha}$ and $B_{\beta}$, and so $\phi$ 
extends uniquely to $B_{\alpha\beta}$.
This extension, also denoted $\phi$, is a 
well-defined amalgam isomorphism from $\B$ to $\C$, and furthermore, 
$\pi_1=\pi_2 \circ \phi$ holds.

In the second case, $B_J\cong C_J\cong G_J$ is isomorphic to a simply connected compact almost simple Lie group of rank one or two. By the universality of the covering $\pi_1 : B_J \rightarrow L_J$, as $B_J$ is simply connected, there exists a unique isomorphism $\psi : B_J \rightarrow C_J$ such that $\pi_1 = \pi_2 \circ \psi$. 
$$\xymatrix{
C_J \ar[dr]_{\pi_2} & B_J \ar[d]^{\pi_1} \ar[l]_{\psi} \\
& L_J
}$$
Consider a 
mapping $\alpha$ from $L_J$ to $L_J$ defined as follows:  For $u\in 
L_J$, let $\alpha(u)=(\pi_2 \circ \psi \circ \pi_1^{-1})(u)$.  Notice that $\alpha$ is a
well-defined automorphism of $L_J$, because the cosets of the kernel of $\pi_1$ are mapped by $\psi$ to cosets of 
the kernel of $\pi_2$. Every automorphism of $L_J$
lifts to a unique automorphism of $C_J$. Indeed, both $L_J$ and $C_J$ are perfect by a corollary of Got\^o's Commutator Theorem (see Corollary 6.56 of \cite{Hofmann/Morris:1998}) and, by Theorem 2.1 of \cite{Sah:1986}, the group $C_J$, which is isomorphic to $\SU_2(\CC) \cong \Spin_3(\RR) \cong \U_1(\HH)$ or to $\SU_3(\CC)$ or to $\Spin_5(\RR) \cong \U_2(\HH)$, is the universal perfect central extension of $L_J$, cf.\ \cite{Huppert:1967} or \cite{Schur:1904}, \cite{Schur:1907}. Alternatively, one can argue as follows:
 Every 
automorphism of $L_J$ is continuous by Corollary 6.56 of \cite{Hofmann/Morris:1998} and van der Waerden's Continuity Theorem (cf.\ Theorem 5.64 of \cite{Hofmann/Morris:1998}), which lifts to a unique continuous automorphism of $C_J$ by \cite{Kawada:1940}, see also \cite{Hofmann:1962}. Finally, this lift in fact is the unique abstract lift of $\alpha$, as any automorphism of $C_J$ again is continuous.

Thus, there is a unique automorphism $\beta$ of $C_J$ such that 
$\pi_2 \circ \beta=\alpha \circ \pi_2$.  Define $\theta : B_J \rightarrow C_J : \theta(b) = (\beta^{-1} \circ \psi)(b)$.  First of all, by definition we have 
$\pi_1|_{B_J}=\pi_2 \circ \theta$, as 
\begin{eqnarray*}
\pi_2 \circ \theta & = & \pi_2 \circ \beta^{-1} \circ \psi \\
& = & \alpha^{-1} \circ \pi_2 \circ \psi \\
& = & \pi_1|_{B_J} \circ \psi^{-1} \circ \pi_2^{-1}|_{L_J} \circ \pi_2 \circ \psi \\
& = & \pi_1|_{B_J}.
\end{eqnarray*}
 Second, for every $J'\subset J$ we have 
that $\theta^{-1} \circ \phi_{|B_{J'}}$ is a lifting to $B_{J'}$ of the 
identity automorphism of $L_{J'}$ and, by the above, it is the identity. For
$\theta^{-1} \circ \phi_{|B_{J'}} = \psi^{-1} \circ \beta \circ \phi_{|B_{J'}}$
and, the following considered on $B_{J'}/\ker({\pi_1}_{|B_{J'}})$,
\begin{eqnarray*}
\psi^{-1} \circ {\pi_2}^{-1}_{|C_{J'}} \circ \alpha \circ \pi_2 \circ \phi_{|B_{J'}} & = & \psi^{-1} \circ {\pi_2}^{-1}_{|C_{J'}} \circ \pi_2 \circ \psi \circ {\pi_1}^{-1}_{|B_{J'}} \circ \pi_2 \circ \phi_{|B_{J'}} \\
& = & {\pi_1}^{-1}_{|B_{J'}} \circ \pi_2 \circ \phi_{|B_{J'}} \\
& = & \id . 
\end{eqnarray*}
This 
shows that $\phi$ and $\theta$ agree on every subgroup $B_{J'}$, which 
allows us to extend $\phi$ to the entire $\hat\B$ by defining it on 
$B_J$ as $\theta$.
Finally, if $\A$ is (strongly) noncollapsing, 
so is its unambiguous covering $\hat\A$, finishing the proof. 
\end{proof}

\section{Uniqueness of unambiguous amalgams} \label{n=3}

Let $\A=(L_{I \backslash \lb i,j \rb})_{(i,j)\in I}$ be an unambiguous strongly noncollapsing irreducible Phan 
amalgam of rank at least two.  We will establish the 
uniqueness of the respective amalgams $\A$ up to isomorphism in a series of lemmas.  The amalgams of rank two are unique by definition. 

\subsection*{Rank three}

Assume the rank of $\mc{A}$ to be three. Since $\A$ is 
unambiguous, each subgroup $L_{I \backslash \{i\}}$ coincides with $L_{I \backslash\{i,j\}} \cap L_{I \backslash\{i,k\}}$ for $\{ i, j, k \} = \{ 1, 2, 3 \}$. 
We want to prove the uniqueness of the amalgam $\mc{A} = (L_{I \backslash \{ i, j\}})_{i, j \in \lb 1, 2, 3 \rb}$. 

\medskip
For $A_3$, i.e., for the diagram $\quad\node_{L_{I \backslash \lb 1 \rb}}\stroke{}\node_{L_{I \backslash \lb 2 \rb}}\stroke{}\node_{L_{I \backslash \lb 3 \rb}}\quad$, recall the isomorphisms 
\begin{eqnarray*}
L_{I \backslash \{2,3\}} & \cong & \SU_3(\CC), \\
L_{I \backslash \{1,3\}} & \cong & \SU_2(\CC) \times \SU_2(\CC), \\
L_{I \backslash \{1,2\}} & \cong & \SU_3(\CC), \\ 
L_{I \backslash \{3\}} = L_{I \backslash \{2,3\}} \cap L_{I \backslash \{1,3\}} & \cong & \SU_2(\CC), \\
L_{I \backslash \{2\}} = L_{I \backslash \{2,3\}} \cap L_{I \backslash \{1,2\}} & \cong & \SU_2(\CC), \\
L_{I \backslash \{1\}} = L_{I \backslash \{1,3\}} \cap L_{I \backslash \{1,2\}} & \cong & \SU_2(\CC). 
\end{eqnarray*}
For $B_3$, i.e., for the diagram $\quad\node_{L_{I \backslash \lb 1 \rb}}\stroke{}\node_{L_{I \backslash \lb 2 \rb}}\dstroke{>}\node_{L_{I \backslash \lb 3 \rb}}\quad$, recall the isomorphisms 
\begin{eqnarray*}
L_{I \backslash \{2,3\}} & \cong & \Spin_5(\RR), \\
L_{I \backslash \{1,3\}} & \cong & \SU_2(\CC) \times \Spin_3(\RR), \\
L_{I \backslash \{1,2\}} & \cong & \SU_3(\CC), \\
L_{I \backslash \{3\}} = L_{I \backslash \{2,3\}} \cap L_{I \backslash \{1,3\}} & \cong & \Spin_3(\RR), \\
L_{I \backslash \{2\}} = L_{I \backslash \{2,3\}} \cap L_{I \backslash \{1,2\}} & \cong & \SU_2(\CC), \\
L_{I \backslash \{1\}} = L_{I \backslash \{1,3\}} \cap L_{I \backslash \{1,2\}} & \cong & \SU_2(\CC).
\end{eqnarray*}
For $C_3$, i.e., for the diagram $\quad\node_{L_{I \backslash \lb 1 \rb}}\stroke{}\node_{L_{I \backslash \lb 2 \rb}}\dstroke{<}\node_{L_{I \backslash \lb 3 \rb}}\quad$, recall the isomorphisms 
\begin{eqnarray*}
L_{I \backslash \{2,3\}} & \cong & \U_2(\HH), \\
L_{I \backslash \{1,3\}} & \cong & \SU_2(\CC) \times \U_1(\HH), \\
L_{I \backslash \{1,2\}} & \cong & \SU_3(\CC), \\
L_{I \backslash \{3\}} = L_{I \backslash \{2,3\}} \cap L_{I \backslash \{1,3\}} & \cong & \U_1(\HH), \\
L_{I \backslash \{2\}} = L_{I \backslash \{2,3\}} \cap L_{I \backslash \{1,2\}} & \cong & \SU_2(\CC), \\
L_{I \backslash \{1\}} = L_{I \backslash \{1,3\}} \cap L_{I \backslash \{1,2\}} & \cong & \SU_2(\CC).
\end{eqnarray*}
Assume there exists another amalgam $\mc{A}' = (L'_{I \backslash \{ i, j\}})_{i, j \in \lb 1, 2, 3 \rb}$. According to Goldschmidt's Lemma (Lemma \ref{Goldschmidt}) the amalgams $\mc{B} = (L_{I \backslash \lb 2,3 \rb}, L_{I \backslash \lb 1,2 \rb}, L_{I \backslash \lb 2 \rb})$ and $\mc{B}' = (L'_{I \backslash \lb 2,3 \rb}, L'_{I \backslash \lb 1,2 \rb}, L'_{I \backslash \lb 2 \rb})$ are isomorphic via some amalgam isomorphism $\psi$, because every automorphism of the group $L_{I \backslash \{2\}} \cong \SU_2(\CC)$ is induced by some automorphism of the group $L_{I \backslash \{1,2\}} \cong \SU_3(\CC)$. Indeed, $L_{I \backslash \{2\}}$ is embedded as the stabilizer of a vector of length one of the natural module of $L_{I \backslash \{1,2\}}$. Clearly, $\psi(L_{I \backslash  \{2\}}) = \psi(L_{I \backslash  \{2,3\}} \cap L_{I \backslash \{1,2\}}) = L'_{I \backslash  \{2,3\}} \cap L'_{I \backslash  \{1,2\}} = L'_{I \backslash  \{2\}}$. The groups $L_{I \backslash \{1\}}$ and $L_{I \backslash \{2\}}$ form a standard pair in $L_{I \backslash \{1,2\}}$, and hence $\psi(L_{I \backslash \{1\}})$ and $L'_{I \backslash \{2\}}=\psi(L_{I \backslash \{2\}})$ form a standard pair in $L'_{I \backslash \{1,2\}}=\psi(L_{I \backslash \{1,2\}})$. Certainly also $L'_{I \backslash \{1\}}$ and $L'_{I \backslash \{2\}}$ form a standard pair in $L'_{I \backslash \{1,2\}}$. Therefore, by Lemma \ref{7.1.1.a}, there exists an automorphism of $L'_{I \backslash \{1,2\}}$ that maps $\psi(L_{I \backslash \{1\}})$ onto $L'_{I \backslash \{1\}}$ and that normalizes $L'_{I \backslash \{2\}}$. Thus, we can assume $\psi(L_{I \backslash \{1\}})=L'_{I \backslash \{1\}}$.

Before we can continue we have to study the amalgam $\mc{A}$ a bit more carefully. 
Define $$D_1 = N_{L_{I \backslash \{1\}}}(L_{I \backslash \{2\}}) \quad \mbox{ and } \quad D_3 = N_{L_{I \backslash \{3\}}}(L_{I \backslash \{2\}})$$ where the groups $L_{I \backslash \{2\}}$, $L_{I \backslash \{1\}}$ are considered as subgroups of $L_{I \backslash \{1,2\}}$ and the groups $L_{I \backslash \{3\}}$, $L_{I \backslash \{2\}}$ are considered as subgroups of $L_{I \backslash \{2.3\}}$. Since $L_{I \backslash \{2\}}$ and $L_{I \backslash \{1\}}$ form a standard pair in $L_{I \backslash \{1,2\}}$, it follows that $D_1$ is a maximal torus in $L_{I \backslash \{1\}} \cong \SU_2(\CC)$. Similarly, $D_3$ is a maximal torus in $L_{I \backslash \{3\}}$. We also define $$D^1_2 = N_{L_{I \backslash \{2\}}}(L_{I \backslash \{1\}}) \quad \mbox{ and } \quad D^3_2 = N_{L_{I \backslash \{2\}}}(L_{I \backslash \{3\}}).$$ Again, these are two maximal tori in $L_{I \backslash \{2\}} \cong \SU_2(\CC)$. The following lemma gives us an extra condition on $\mc{A}$ that holds because $\mc{A}$ is strongly noncollapsing.

\begin{lemma} \label{non}\label{non-collapsing}
$D^1_2 = D^3_2$.
\end{lemma}

\begin{proof}
Let $G$ be a nontrivial completion of $\mc{A}$ and let $\pi$ be the corresponding map from $\mc{A}$ to $G$. Since $\mc{A}$ is assumed to be strongly noncollapsing, we may assume that $\pi$ is injective on every $L_{I \backslash \{i\}}$. Observe that $D^i_2 = C_{L_{\{1,3\}}}(D_i)$ for $i = 1, 3$. Thus, $\pi(D_2^i) = C_{\pi(L_{I \backslash \{2\}})}(\pi(D_i))$. Since $D_1$ and $D_3$ commute elementwise in $L_{I \backslash \{1,3\}}$, we have that $\pi(D_1)$ and $\pi(D_3)$ commute elementwise as well. Since $L_{I \backslash \{2\}}$ is invariant under $D_1=N_{L_{I \backslash \{1\}}}(L_{I \backslash \{2\}})$ (in $L_{I \backslash \{1,2\}}$) and since $\pi$ is injective on $L_{I \backslash \{2\}}$, it follows that $D^3_2=C_{L_{I \backslash \{2\}}}(D_3)$ is invariant under $D_1$ (again as subgroups of $L_{I \backslash \{1,2\}}$) and $\pi(D_2^3) = C_{\pi(L_{I \backslash \{2\}})}(\pi(D_3))$ is invariant under $\pi(D_1)$. Here, injectivity of $\pi$ is needed for the following argument. $D_1$ and $D_3$ commute as subgroups of $L_{I \backslash \{1,3\}}$. The group $L_{I \backslash \{2\}}$ is invariant under $D_1$ as a subgroup of $L_{I \backslash \{1,2\}}$. Since $L_{I \backslash \{1,3\}}$ and $L_{I \backslash \{1,2\}}$ are not contained in a common group of the amalgam $\mc{A}$, we cannot conclude that $D_1$ leaves $D^3_2$ invariant. However, in $G$, since $L_{I \backslash \{2\}}$, $D_1$, $D_3$, $D^3_2$ are embedded via $\pi$, we can draw that conclusion. 

But now the maximal torus $D_1$ of $L_{I \backslash \{1\}} \cong \SU_2(\CC)$ leaves invariant the maximal tori $D^1_2$ and $D^3_2$ of $L_{I \backslash \{2\}} \cong \SU_2(\CC)$. Analysis of the group $L_{I \backslash \{1,2\}} \cong \SU_3(\CC)$ shows that $D_2^1 = D^3_2$. 
\end{proof}

In view of this lemma we can use the notation $$D_2=D_2^1=D_2^3.$$
Since $N_{L_{I \backslash \{2\}}}(L_{I \backslash \{1\}}) = D^1_2 = D_2 = D^3_2 = N_{L_{I \backslash \{2\}}}(L_{I \backslash \{3\}})$, the considerations made before Lemma \ref{non-collapsing} imply $\psi(D_2)=D_2'$.
Let $d$ be a nontrivial element of $D_2'$ of order distinct from two. Denote by $W$ the natural three-dimensional module of $L'_{I \backslash \{1,2\}}$, and recall that $L'_{I \backslash \{2\}}$ and $L'_{I \backslash \{3\}}$ form a standard pair of $L'_{I \backslash \{2,3\}}$. As $D_2' \leq L'_{I \backslash \{2\}}$, the group $D_2'$ fixes a non-isotropic vector $u$ of length one of $W$ fixed by $L'_{I \backslash \{2\}}$. Since $D_2'$ normalizes $L'_{I \backslash \{3\}}$, it also stabilizes $\gen{v}$, where $v$ is a non-isotropic vector of length one of $W$ fixed by $L'_{I \backslash \{3\}}$. Moreover, since $L'_{I \backslash \{2\}}$ and $L'_{I \backslash \{3\}}$ form a standard pair, $u$ is perpendicular to $v$ in $W$. Let $\gen{w}$ be the one-dimensional subspace of $W$ that is perpendicular to both $u$ and $v$ and assume $w$ has length one. Then $u$, $v$, $w$ is an orthonormal basis of $W$, and $d$ acts diagonally with respect to that basis via $\diag(1,a,a^{-1})$. Since the order of $d$ is distinct from two, we have $a \neq a^{-1}$, so the one-dimensional subspaces of $W$ stabilized by $d$ are precisely $\gen{u}$, $\gen{v}$, $\gen{w}$. It follows, since $D_2' = \psi(D_2) = N_{\psi(L_{I \backslash \{2\}})}(\psi(L_{I \backslash \{3\}})) = N_{L'_{I \backslash \{2\}}}(\psi(L_{I \backslash \{3\}}))$, that $\psi(L_{I \backslash \{3\}})$ is the stabilizer of either $v$ or $w$. 

In the former case we have $\psi(L_{I \backslash \{3\}}) = L'_{I \backslash \{3\}}$, and we have proved $\mc{A} \cong \mc{A}'$, since $L_{I \backslash \{1,3\}} = L_{I \backslash \{3\}} \times L_{I \backslash \{1\}}$ and $L'_{I \backslash \{1,3\}} = L'_{I \backslash \{3\}} \times L'_{I \backslash \{1\}}$.

In the latter case consider the element $g$ of $L'_{I \backslash \{2\}}$ whose matrix with respect to the orthonormal basis $u$, $v$, $w$ has the form
$$
\begin{pmatrix}
1 & 0 & 0 \\
0 & 0 & -1 \\
0 & 1 & 0
\end{pmatrix} . $$
Conjugation with $g$ induces the action of the contragredient automorphism on $L'_{I \backslash \{2\}}$. By the defining relation $$A^{-1} = \bar A^{T}$$ of unitary matrices the action of the contragredient automorphism of $L'_{I \backslash \{2\}}$ coincides with the field involution. Therefore, we can define an automorphism $\alpha$ of $\mc{B}'$ that acts trivially on $L'_{I \backslash \{1,2\}}$ and as the composition of the field automorphism and conjugation by $g$ on $L'_{I \backslash \{2,3\}}$, since by the above this automorphism acts trivially on $L'_{I \backslash \{2\}} = L'_{I \backslash \{2,3\}} \cap L'_{I \backslash \{1,2\}}$. Moreover, $\alpha$ interchanges $\gen{v}$ and $\gen{w}$, so it maps $\psi(L_{I \backslash \{3\}})$ onto $L'_{I \backslash \{3\}}$.

\medskip
We have proved the following.

\begin{proposition} \label{rank 3 case}
Let $\mc{A}$ be a strongly noncollapsing unambiguous irreducible Phan amalgam of rank three. Then $\A$ is unique up to isomorphism, i.e., $\A$ is isomorphic to a standard Phan amalgam. \pend
\end{proposition}

\subsection*{Rank at least four}

Let 
$\A = (L_{I \backslash \lb i, j \rb})_{1 \leq i < j \leq n}$ be a strongly noncollapsing unambiguous irreducible Phan amalgam of rank at least four.
We complete the proof of the uniqueness of $\A$ by induction, the case of rank three from Proposition \ref{rank 3 case} being the basis of induction. 

\begin{lemma} \label{extension}\label{101}
Let $n \geq 4$ and let 
$\A$ be a strongly noncollapsing unambiguous irreducible Phan amalgam of rank $n$.
Then there exists a unique amalgam 
\begin{eqnarray*}
\B_\A & = & \A\cup H_1\cup H_2 \\ 
& \mbox{ with } & \\
H_1 & = & \gen{L_{I \backslash \lb i, j \rb} \mid 1\le i<j\le n-1}  \mbox{ and } \\
H_2 & = & \gen{L_{I \backslash \lb i, j \rb} \mid 2\le i<j\le n}  .
\end{eqnarray*}
The group $H_1$ is isomorphic to $\SU_n(\mathbb{C})$ unless the case of the Dynkin diagram $F_4$, where $H_1$ is isomorphic to $\Spin_7(\mathbb{R})$, while the group $H_2$ is isomorphic to
\begin{eqnarray*}
\SU_{n}(\CC) & & \mbox{ for the diagram $A_n$}, \\
\Spin_{2n-1}(\RR) & & \mbox{ for the diagram $B_n$}, \\
\U_{n-1}(\HH) & & \mbox{ for the diagram $C_n$}, \\
\Spin_{2n-2}(\RR) & & \mbox{ for the diagram $D_n$}, \\
\Spin_{10}(\RR) & & \mbox{ for the diagram $E_6$}, \\
\Spin_{12}(\RR) & & \mbox{ for the diagram $E_7$}, \\
\Spin_{14}(\RR) & & \mbox{ for the diagram $E_8$}, \\
\U_{3}(\HH) & & \mbox{ for the diagram $F_4$}.
\end{eqnarray*}
\end{lemma}

\begin{proof} Let 
\begin{eqnarray*}
\B_1 & := & (L_{I \backslash \lb i, j \rb})_{1\le i<j\le n-1}, \\ 
\B_2 & := & (L_{I \backslash \lb i, j \rb})_{2\le i<j\le n}, \quad \mbox{ and } \\ 
\C & := & \B_1\cap\B_2 .
\end{eqnarray*}
By the inductive assumption, 
both $\B_1$ and $\B_2$ are isomorphic to some standard Phan amalgam
 and hence there exist faithful completions $\pi_i:\B_i \rightarrow H_i$ where the isomorphism types of
$H_1$ and $H_2$ are given as in the hypothesis. We want to glue $H_1$ and $H_2$ to the amalgam $\A$ via 
$\pi_1$ and $\pi_2$. Let $K_i := \gen{\pi_i(\mc{C})}$. Since, again by the inductive assumption, the amalgam $\mc{C}$ is isomorphic to a standard Phan amalgam, we have $K_i \cong \SU_{n-1}(\CC)$ or, in case of the diagram $F_4$, we have $K_i \cong \Spin_{5}(\RR) \cong \U_2(\HH)$. By Proposition \ref{characteristic} the group $K_i$ is a characteristic completion of the amalgam $\mc{C}$, so there exists an isomorphism 
$\phi:K_1\rightarrow K_2$ that takes $\pi_1(\C)$ to $\pi_2(\C)$. 
Let $\psi$ be the restriction of $\phi$ to $\pi_1(\C)$.  Applying the Bennett-Shpectorov Lemma (Lemma \ref{64}) with $\phi : K_1 \rightarrow K_2$ and $\psi : \pi_1 (\C) \rightarrow \pi_2(\C)$ as above and $\psi' : \pi_1 (\C) \rightarrow \pi_2(\C)$ with $\psi' = \pi_2 \circ {\pi_1}_{|\mc{C}}^{-1}$, there exists a unique isomorphism $\phi' : K_1 \rightarrow K_2$ such that $\phi'_{|\pi_1(\mc{C})} = \psi'$. Thus, $\phi' \circ {\pi_1}_{|\mc{C}} = \pi_2|_{\mc{C}}$. Identifying $K_1$ with $K_2$ via $\phi'$ we obtain our unique amalgam $\mc{B}$.
\end{proof}

Let us now turn to the uniqueness of the amalgam $\mc{A}$. Suppose we 
have strongly noncollapsing unambiguous irreducible Phan amalgams $\A$ and $\A'$ corresponding to the same diagram.  Extend $\A$ and $\A'$ to 
amalgams $\B_\A=\A\cup H_1\cup H_2$ and $\B'_{\A'}=\A'\cup H'_1\cup H'_2$ as in 
Lemma \ref{extension}. By Goldschmidt's Lemma (Lemma \ref{Goldschmidt}) there 
exists an isomorphism $\phi$ from $H_1\cup H_2$ onto $H'_1\cup H'_2$.  
By the inductive assumption $(L_{I \backslash \lb i, j \rb})_{1<i<j<n}$ is isomorphic to a standard Phan amalgam embedded in $H_1\cap H_2$.  Similarly $(L'_{I \backslash \lb i, j \rb})_{1<i<j<n}$ and $\phi(L_{I \backslash \lb i, j \rb})_{1<i<j<n}$ are isomorphic to standard Phan amalgams embedded in $H'_1\cap H'_2$.  These two amalgams 
correspond to two choices of a maximal torus of $H'_1\cap H'_2$, which are conjugate by Theorem 6.27 of \cite{Hofmann/Morris:1998}.  So, correcting $\phi$ if 
necessary by an inner automorphism of $H'_1\cap H'_2$, we may assume 
that $\phi(L_{I \backslash \lb i \rb }) = L'_{I \backslash \lb i \rb }$ for $ 1 < i < n$ and $\phi(L_{I \backslash \lb i, j \rb})=L'_{I \backslash \lb i, j \rb}$ for $1<i<j<n$. Also, by studying the standard Phan amalgam inside $H_1'$, we have
\begin{eqnarray*}
\phi\left(L_{I \backslash \{1 \}}  \right) & = & \phi\left(C_{H_1}\left(\gen{L_{I \backslash \{3\}}, \ldots, L_{I \backslash \{n-1 \}}}\right)\right) \\
& = & C_{\phi(H_1)}\left(\phi\left(\gen{L_{I \backslash \{3\}}, \ldots, L_{I \backslash \{n-1 \}}}\right)\right) \\ 
& = & C_{H'_1}\left(\gen{L'_{I \backslash \{3\}}, \ldots, L'_{I \backslash \{n-1 \}}}\right)  \\
& = & L'_{I \backslash \{1 \}}.
\end{eqnarray*}
By a similar argument, $\phi(L_{I \backslash \{n \}}) = L'_{I \backslash \{n \}}$. Therefore $\phi$ extends to an isomorphism from $\mc{A}$ to $\mc{A}'$. Indeed, $\phi$ is already defined on all $L_{I \backslash \{i,j \}}$ with $2 \leq i < j \leq n-1$. Also, inside the standard Phan amalgam of $H_1'$ we see that $\phi(L_{I \backslash \{1,i \}}) = L'_{I \backslash \{1,i \}}$ for $i < n$, since $L_{I \backslash \{1,i \}} = \gen{L_{I \backslash \{1 \}},L_{I \backslash \{i \}}}$. Similarly, in the standard Phan amalgam of $H_2'$ we see that $\phi(L_{I \backslash \{i,n \}}) = L'_{I \backslash \{i,n \}}$ for $1 < i$. It remains to realize that $L_{I \backslash \{1,n\}}$ is the direct product of $L_{I \backslash \{1 \}}$ and $L_{I \backslash \{n \}}$, so that $\phi$ extends to an isomorphism of $\mc{A}$ to $\mc{A}'$. 

\medskip Thus we have shown:

\begin{proposition} \label{n>=4}
Let $n \geq 4$, and let 
$\A$ be a strongly noncollapsing unambiguous irreducible Phan amalgam of rank $n$. Then $\mc{A}$ is unique up to isomorphism, i.e., $\mc{A}$ is isomorphic to a standard Phan amalgam. \pend
\end{proposition}

\medskip
\noindent
{\bf Proof of the Main Theorem.}
The weak Phan system of $G$ gives rise to a strongly noncollapsing Phan amalgam $\mc{A}$, which by Proposition \ref{similar} is covered by a unique strongly noncollapsing unambiguous Phan amalgam $\widehat{\mc{A}}$. This strongly noncollapsing unambiguous Phan amalgam $\widehat{\mc{A}}$ is isomorphic to a standard Phan amalgam by Propositions \ref{rank 3 case} and \ref{n>=4} applied to the irreducible components of $\Delta$ of rank at least three and by Definition \ref{arbitraryan} applied to the irreducible components of $\Delta$ of rank at most two. Finally, the first claim follows by Theorem \ref{borovoi}. The second claim follows immediately from the first claim by the classification of irreducible Dynkin diagrams, see \cite{Bourbaki:2002}, and by \cite{Helgason:1978} or by 94.33 of \cite{Salzmann:1995}.
\pend


\vspace{2cm}

\noindent Author's address:

\medskip
\noindent Ralf Gramlich \\
TU Darmstadt \\
FB Mathematik / AG 5 \\
Schlo\ss gartenstra\ss e 7 \\
64289 Darmstadt \\
Germany \\
{\tt gramlich@mathematik.tu-darmstadt.de}

\end{document}